# A note on sensitivity of principal component subspaces and the efficient detection of influential observations in high dimensions

**Luke A. Prendergast**

*La Trobe University*
*Dept of Mathematics and Statistics,*
*La Trobe University, VIC 3086, Australia.*
*e-mail:* luke.prendergast@latrobe.edu.au

**Abstract:** In this paper we introduce an influence measure based on second order expansion of the RV and GCD measures for the comparison between unperturbed and perturbed eigenvectors of a symmetric matrix estimator. Example estimators are considered to highlight how this measure compliments recent influence analysis. Importantly, we also show how a sample based version of this measure can be used to accurately and efficiently detect influential observations in practice.

## Contents

## 1. Introduction

    Principal component analysis (PCA) is a popular tool in multivariate statistics. However, PCA estimates may be highly influenced by certain types of

observations and, as such, it is often important to locate, and perhaps subsequently treat, those observations which may be potentially harmful in practice. Influence analysis of PCA estimators (see, for e.g., [9; 3; 31; 11] and [14]) show that the direction of an observation plays a crucial role in how influential it may be on the eigenvector and eigenvalue estimates. As such, influential observations may or may not be detectable using common distance measures and can therefore be difficult to locate.

Throughout we will consider PCA of a symmetric matrix where the set of all eigenvectors satisfies orthonormality conditions. Often it is the span of a subset of eigenvectors that is of primary interest rather than individual elements. Let $\mathbf{A} = [\mathbf{a}_1, \ldots, \mathbf{a}_K]$ and $\mathbf{B} = [\mathbf{b}_1, \ldots, \mathbf{b}_K]$ be two $p \times K$ matrices where $\|\mathbf{a}_i\| = 1$, $\|\mathbf{b}_i\| = 1$ $(i = 1, \ldots, K)$ and $\mathbf{a}_i^\top \mathbf{a}_j = 0$, $\mathbf{b}_i^\top \mathbf{b}_j = 0$ for all $i \neq j$. Two common measures for the comparison between the column spaces of $\mathbf{A}$ and $\mathbf{B}$ are the RV$(\mathbf{A}, \mathbf{B})$ coefficient ([15; 27]) and the GCD$(\mathbf{A}, \mathbf{B})$ measure [32]. Let $\mathbf{P_A} = \mathbf{A}\mathbf{A}^\top$ and $\mathbf{P_B} = \mathbf{B}\mathbf{B}^\top$ denote the projection matrices onto the column spaces of $\mathbf{A}$ and $\mathbf{B}$ respectively. Then, due to orthonormality of the columns of $\mathbf{A}$ and $\mathbf{B}$ (see, for e.g., [3]) we have

$$\text{RV}(\mathbf{A}, \mathbf{B}) = \text{GCD}(\mathbf{A}, \mathbf{B}) = \frac{1}{K}\text{trace}\left(\mathbf{P_A}\mathbf{P_B}\right). \tag{1}$$

In considering influence on eigenvector subset spans, Bénasséni [3] noted that the RV and GCD measures were insensitive to small perturbations. Bénasséni then introduced a new measure for assessing sensitivity given as

$$\rho_1(\mathbf{A}, \mathbf{B}) = 1 - \frac{1}{K}\sum_{k=1}^{K} \|\mathbf{a}_k - \mathbf{P_B}\mathbf{a}_k\| \tag{2}$$

or, alternatively, $\rho_2(\mathbf{A}, \mathbf{B}) = \rho_1(\mathbf{B}, \mathbf{A})$. It should be noted that it is not necessarily true that $\rho_2(\mathbf{A}, \mathbf{B}) = \rho_1(\mathbf{A}, \mathbf{B})$ so that, unlike the RV and GCD measures, the ordering of the arguments may be important. However, by considering a small adjustment in $\rho_1(\mathbf{A}, \mathbf{B})$ that considers summation of $\|\mathbf{a}_k - \mathbf{P_B}\mathbf{a}_k\|^2$ instead of $\|\mathbf{a}_k - \mathbf{P_B}\mathbf{a}_k\|$, note that

$$1 - \frac{1}{K}\sum_{k=1}^{K} \|\mathbf{a}_k - \mathbf{P_B}\mathbf{a}_k\|^2 = \text{RV}(\mathbf{A}, \mathbf{B}) = \text{GCD}(\mathbf{A}, \mathbf{B})$$

so that there is a strong link between Bénasséni's measure and the RV and GCD measures. The purpose of this paper is then to consider an influence measure based on the RV and GCD measures.

In Section 2 we consider perturbation of the RV and GCD measures and introduce a measure for the influence analysis of eigenvector spans. We then apply this measure to some example estimators in Section 3 to show how it compliments existing influence studies. Sample versions are discussed in Section 4 for the detection of highly influential observations in practice. In Section 5 we show how influential observations may be efficiently detected in practice with respect to a high dimensional data set.

## 2. The effect of perturbation on the RV and GCD measures

Consider the contamination distribution

$$F_{\mathbf{x}}(\varepsilon) = (1-\varepsilon)F_{\boldsymbol{\mu},\boldsymbol{\Sigma}} + \varepsilon\Delta_{\mathbf{x}} \tag{3}$$

where $F_{\boldsymbol{\mu},\boldsymbol{\Sigma}}$ is some distribution with mean $\boldsymbol{\mu}$ and covariance matrix $\boldsymbol{\Sigma}$, $0 < \varepsilon < 1$ and $\Delta_{\mathbf{x}}$ is the Dirac measure putting all of its probability mass at the contaminant $\mathbf{x}$. When convenient we may also utilize $\mathbf{z} = \boldsymbol{\Sigma}^{-1/2}(\mathbf{x} - \boldsymbol{\mu})$ (the standardized contaminant at $F_{\boldsymbol{\mu},\boldsymbol{\Sigma}}$) and the population Mahalanobis distance with respect to $\mathbf{x}$ at $F_{\boldsymbol{\mu},\boldsymbol{\Sigma}}$ given as

$$\mathrm{MD}_{\boldsymbol{\mu},\boldsymbol{\Sigma}}(\mathbf{x}) = \sqrt{(\mathbf{x}-\boldsymbol{\mu})^{\top}\boldsymbol{\Sigma}^{-1}(\mathbf{x}-\boldsymbol{\mu})} = \|\mathbf{z}\|.$$

Let $W$ denote a $p \times p$ statistical functional where, at an arbitrary distribution $G$ for which it exists, $W(G)$ is symmetric. With respect to $F_{\mathbf{x}}(\varepsilon)$ defined in (3), we are interested in perturbation of the form

$$W(F_{\mathbf{x}}(\varepsilon)) = W(F_{\boldsymbol{\mu},\boldsymbol{\Sigma}}) + \varepsilon\mathbf{W}_1 + \varepsilon^2\mathbf{W}_2 + O(\varepsilon^3) \tag{4}$$

where $\mathbf{W}_1$ and $\mathbf{W}_2$ are $p \times p$ symmetric matrices independent of $\varepsilon$.

Under perturbation of the form given in (4), the first order coefficient of $\varepsilon$ is the influence function ([16; 17]) for $W$ at $F_{\boldsymbol{\mu},\boldsymbol{\Sigma}}$ denoted $\mathrm{IF}(W, F_{\boldsymbol{\mu},\boldsymbol{\Sigma}}; \mathbf{x})$. The influence function is a useful tool for understanding the sensitivity of estimators to small perturbations. For example, if $\mathrm{IF}(W, F_{\boldsymbol{\mu},\boldsymbol{\Sigma}}; \mathbf{x}) = 0$ then, for small $\varepsilon$, $W(F_{\mathbf{x}}(\varepsilon)) \approx W(F_{\boldsymbol{\mu},\boldsymbol{\Sigma}})$ so that small perturbations with respect to $\mathbf{x}$ have little influence on the estimator. On the other hand, if $\mathrm{IF}(W, F_{\boldsymbol{\mu},\boldsymbol{\Sigma}}; \mathbf{x})$ is large then perturbation with respect to $\mathbf{x}$ is highly influential since it effects a large change on the estimator.

Let $|\kappa_1| \geq |\kappa_2| \geq \cdots \geq |\kappa_p|$ be the eigenvalues of $W(F_{\boldsymbol{\mu},\boldsymbol{\Sigma}})$ and let $\boldsymbol{\nu}_1, \ldots, \boldsymbol{\nu}_p$ denote the corresponding eigenvectors. We are interested in the effect that perturbation has on the span of a subset of the eigenvectors. Let $S \subset \{1, \ldots, p\}$ and let $\mathbf{P}_S = \sum_{j \in S} \boldsymbol{\nu}_j \boldsymbol{\nu}_j^{\top}$ denote the projection matrix onto the subspace spanned by the elements of $\{\boldsymbol{\nu}_j : j \in S\}$. Similarly, let $\mathbf{P}_S(\varepsilon) = \sum_{j \in S} \boldsymbol{\nu}_j(\varepsilon)\boldsymbol{\nu}_j(\varepsilon)^{\top}$ denote the perturbed equivalent at $F_{\mathbf{x}}(\varepsilon)$ where $\boldsymbol{\nu}_1(\varepsilon), \ldots, \boldsymbol{\nu}_p(\varepsilon)$ are the eigenvectors corresponding to the ordered absolute eigenvalues of $W(F_{\mathbf{x}}(\varepsilon))$. Typically, $S$ will be chosen to be $\{1, \ldots, K\}$ such that the span of the eigenvectors corresponding to the $K$ largest absolute eigenvalues is of interest. For example, in PCA corresponding to covariance matrices, principal components with corresponding large eigenvalues are retained as they can account for most of the total population variance. Let $S'$ denote the compliment of $S$. The following condition will also be used.

**Condition 1.** *For each $j \in S$ and $r \in S'$, $\kappa_j \neq \kappa_r$.*

We now look at a measure based on the expansion of a function of the RV and GCD measures based on the above condition. The proof is given in the Appendix.

**Theorem 1.** *Consider $S$, $\mathbf{P}_S$, $\mathbf{P}_S(\varepsilon)$, $\kappa_i$ $(i = 1, \ldots, p)$ and $\boldsymbol{\nu}_i$ $(i = 1, \ldots, p)$ defined previously and let*

$$\rho_S(W, F_{\boldsymbol{\mu}, \boldsymbol{\Sigma}}; \mathbf{x}) = \frac{1}{\varepsilon^2} \left[ 1 - \frac{1}{K} \text{trace} \left\{ \mathbf{P}_S \mathbf{P}_S(\varepsilon) \right\} \right].$$

*Then, under the perturbation form given in* (4) *and Condition* 1,

$$\rho_S(W, F_{\boldsymbol{\mu}, \boldsymbol{\Sigma}}; \mathbf{x}) = \frac{1}{K} \sum_{j \in S} \sum_{r \in S'} \frac{\{\boldsymbol{\nu}_j^{\top} \text{IF}(W, F_{\boldsymbol{\mu}, \boldsymbol{\Sigma}}; \mathbf{x}) \boldsymbol{\nu}_r\}^2}{(\kappa_j - \kappa_r)^2} + O(\varepsilon)$$

*where $K$ is the number of elements in the set $S$ and $\text{IF}(W, F_{\boldsymbol{\mu}, \boldsymbol{\Sigma}}; \mathbf{x})$ is the influence function for $W$ at $F_{\boldsymbol{\mu}, \boldsymbol{\Sigma}}$.*

From Theorem 1, if perturbation has resulted in no difference between the span of the non-perturbed and perturbed eigenvectors then $\rho_S(W, F_{\boldsymbol{\mu}, \boldsymbol{\Sigma}}; \mathbf{x}) = 0$ since $\text{trace} \left\{ \mathbf{P}_S \mathbf{P}_S(\varepsilon) \right\} = K$. However, as the distance between the spans increases according to

$$\frac{1}{K} \text{trace} \left\{ \mathbf{P}_S \mathbf{P}_S(\varepsilon) \right\}$$

(i.e, as the RV and GCD measure for the comparison between the spans approaches zero) then $\rho_S(W, F_{\boldsymbol{\mu}, \boldsymbol{\Sigma}}; \mathbf{x})$ increases. Also, for small $\varepsilon$ we have that $\rho_S(W, F_{\boldsymbol{\mu}, \boldsymbol{\Sigma}}; \mathbf{x})$ is approximately equal to the second order coefficient to $\varepsilon$ in the expansion of $[1 - \text{trace} \left\{ \mathbf{P}_S \mathbf{P}_S(\varepsilon) \right\} / K]$. Hence, the following influence measure will be utilized throughout,

$$\tilde{\rho}_S(W, F_{\boldsymbol{\mu}, \boldsymbol{\Sigma}}; \mathbf{x}) = \lim_{\varepsilon \downarrow 0} \rho_S(W, F_{\boldsymbol{\mu}, \boldsymbol{\Sigma}}; \mathbf{x}) = \frac{1}{K} \sum_{j \in S} \sum_{r \in S'} \frac{\{\boldsymbol{\nu}_j^{\top} \text{IF}(W, F_{\boldsymbol{\mu}, \boldsymbol{\Sigma}}; \mathbf{x}) \boldsymbol{\nu}_r\}^2}{(\kappa_j - \kappa_r)^2}.$$

(5)

In the case of $K = 1$ such that $j \in \{1, \ldots, p\}$,

$$\tilde{\rho}_S(W, F_{\boldsymbol{\mu}, \boldsymbol{\Sigma}}; \mathbf{x}) = \|\text{IF}(\nu_j, F_{\boldsymbol{\mu}, \boldsymbol{\Sigma}}; \mathbf{x})\|^2$$

where $\nu_j$ is the functional for the $j$th eigenvector estimator and $\text{IF}(\nu_j, F_{\boldsymbol{\mu}, \boldsymbol{\Sigma}}; \mathbf{x})$ is the associated influence function. Hence, the measure may be used to assess influence on individual components as well as the span of a subset of components.

The $\tilde{\rho}_S(W, F_{\boldsymbol{\mu}, \boldsymbol{\Sigma}}; \mathbf{x})$ measure provides a convenient means to understand the sensitivity of eigenvector estimators in the presence of non-unique eigenvalues. $\text{IF}(\nu_j, F_{\boldsymbol{\mu}, \boldsymbol{\Sigma}}; \mathbf{x})$ is only known for the case of unique $\kappa_j$ (see, for e.g. [9] or [11]) and, as such, is not useful in all situations. For example, let $\mathcal{J} \subset \{1, \ldots, p\}$ such that $\{\kappa_i = \kappa_j : i, j \in \mathcal{J}, i \neq j\}$ then the solution to $\text{IF}(\nu_j, F_{\boldsymbol{\mu}, \boldsymbol{\Sigma}}; \mathbf{x})$ is not known for $j \in \mathcal{J}$. However, from (5), $\tilde{\rho}_S(W, F_{\boldsymbol{\mu}, \boldsymbol{\Sigma}}; \mathbf{x})$ is known provided $\mathcal{J} \subset S$ or $\mathcal{J} \subset S'$ and Condition 1 holds.

## 3. Example estimators

### 3.1. Covariance and correlation matrix estimators

Let $C_0$ denote the functional for the classical covariance matrix estimator where, at $F_{\boldsymbol{\mu}, \boldsymbol{\Sigma}}$, $C_0(F_{\boldsymbol{\mu}, \boldsymbol{\Sigma}}) = \boldsymbol{\Sigma}$ with eigenvalues $\lambda_1 \geq \cdots \geq \lambda_p$ and correspond-ing eigenvectors $\boldsymbol{\eta}_1, \ldots, \boldsymbol{\eta}_p$. The influence function for this estimator (see, for e.g. [9]) is $\text{IF}(C_0, F_{\boldsymbol{\mu}, \boldsymbol{\Sigma}}; \mathbf{x}) = (\mathbf{x} - \boldsymbol{\mu})(\mathbf{x} - \boldsymbol{\mu})^\top - \boldsymbol{\Sigma}$ so that, from (5),

$$\tilde{\rho}_S(C_0, F_{\boldsymbol{\mu}, \boldsymbol{\Sigma}}; \mathbf{x}) = \frac{1}{K} \sum_{j \in S} \sum_{r \in S'} \frac{y_j^2 y_r^2}{(\lambda_j - \lambda_r)^2}. \tag{6}$$

where $y_j = \boldsymbol{\eta}_j^\top (\mathbf{x} - \boldsymbol{\mu})$.

**Remark 1.** *Bénasséni's coefficient, as shown in (2), was introduced in the classical covariance matrix setting. This measure was based on the average sine of the angle between each of the perturbed eigenvectors and their projection onto the non-perturbed subspace. Using our notation the influence measure associated with this coefficient is (see [3])*

$$\frac{1}{K} \sum_{j \in S} \left\{ \sum_{r \in S'} \frac{y_j^2 y_r^2}{(\lambda_j - \lambda_r)^2} \right\}^{1/2}$$

*which contains similar sensitivity information as the $\tilde{\rho}_S(C_0, F_{\boldsymbol{\mu}, \boldsymbol{\Sigma}}; \mathbf{x})$ measure.*

Similarly, let $R_0$ denote the functional for the classical correlation matrix estimator where, at $F_{\boldsymbol{\mu}, \boldsymbol{\Sigma}}$, $R_0(F_{\boldsymbol{\mu}, \boldsymbol{\Sigma}}) = \boldsymbol{\Gamma}$ and let $\alpha_1, \ldots, \alpha_p$ and $\boldsymbol{\gamma}_1, \ldots, \boldsymbol{\gamma}_p$ denote the eigenvalues and associated eigenvectors of $\boldsymbol{\Gamma}$. For $x_i$ denoting the $i$th element of $\mathbf{x}$, $\mu_i$ denoting the $i$th element of $\boldsymbol{\mu}$ and $\sigma_{ii}$ denoting the $i$th diagonal element of $\boldsymbol{\Sigma}$, let $\tilde{\mathbf{z}}^\top = [\tilde{z}_1, \ldots, \tilde{z}_p] = [(x_1 - \mu_1)/\sigma_{11}, \ldots, (x_p - \mu_p)/\sigma_{pp}]$ and $\mathbf{D} = \text{diag}(\tilde{z}_1^2, \ldots, \tilde{z}_p^2)$. The influence function for this estimator is, from [13], $\text{IF}(R_0, F_{\boldsymbol{\mu}, \boldsymbol{\Sigma}}; \mathbf{x}) = \tilde{\mathbf{z}}\tilde{\mathbf{z}}^\top - (\mathbf{D}\boldsymbol{\Gamma} + \boldsymbol{\Gamma}\mathbf{D})/2$ so that, from (5),

$$\tilde{\rho}_S(R_0, F_{\boldsymbol{\mu}, \boldsymbol{\Sigma}}; \mathbf{x}) = \frac{1}{K} \sum_{j \in S} \sum_{r \in S'} \frac{1}{(\alpha_j - \alpha_r)^2} \left\{ u_j u_r - \frac{1}{2} (\alpha_j + \alpha_r) \boldsymbol{\gamma}_j^\top \mathbf{D} \boldsymbol{\gamma}_r \right\}. \tag{7}$$

We will now provide an example of the form of the influence measure for eigenvector subspaces of covariance estimators. We will not only consider the measure for the classical case as shown in (6), but also with respect to two robust estimators of the covariance matrix; namely the one-step re-weighted Minimum Covariance Determinant (RMCD) estimator which includes an initial MCD estimation step [28] followed by a subsequent re-weighting [21] and the $S$-estimator ([29; 30; 12]). For simplicity and to satisfy Fisher consistency at the non-contaminated model we will assume $F_{\boldsymbol{\mu}, \boldsymbol{\Sigma}}$ is multivariate normal. For the RMCD estimator we choose the breakdown point for the initial MCD estimator to be $\alpha = 0.5$ followed by the retention of 97.5% mass that satisfies $\text{MD}^2_{\boldsymbol{\mu}^*, \boldsymbol{\Sigma}^*}(\mathbf{x}) \le q_{.975}$ where $P(\chi_p^2 \le q_{.975}) = 0.975$ and $\boldsymbol{\mu}^*$ and $\boldsymbol{\Sigma}^*$ are the MCD mean vector and covariance matrix. For the $S$-estimator we use the minimizer function associated with Tukey's biweight function (see, for e.g., Example 2.2 of [20]) and $\alpha = 0.5$. Associated influence functions for these estimators that are used in the following example can be found in [10] and [20].

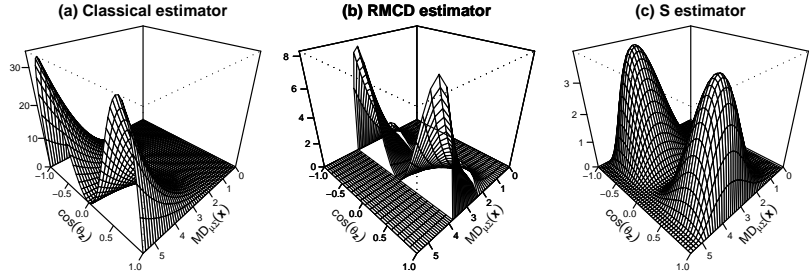

FIG 1. *Subspace sensitivity plots for Example 1 with $K = 3$ for (a) the classical estimator (b) RMCD estimator with $\alpha = 0.5$ and $\delta = 0.025$ and (c) the S-estimator.*

**Example 1.** Let $\lambda_1 = \cdots = \lambda_K$ and $\lambda_{K+1} = \cdots = \lambda_p$ with $\lambda_p < \lambda_1$, $\boldsymbol{\Sigma} = \operatorname{diag}(\lambda_1, \ldots, \lambda_p)$, $\boldsymbol{\mu} = \mathbf{0}$ and $S = \{1, \ldots, K\}$. From (6), and since each $y_i = \lambda_i^{1/2} \boldsymbol{\eta}_i^\top \mathbf{z}$ for $\mathbf{z} = \boldsymbol{\Sigma}^{-1/2} \mathbf{x}$, then $\tilde{\rho}_S(C_0, F_{\boldsymbol{\mu}, \boldsymbol{\Sigma}}; \mathbf{x})$ is equal to

$$\lambda_1 \lambda_p \{K(\lambda_1 - \lambda_p)^2\}^{-1} \operatorname{trace}(\mathbf{P}_S \mathbf{z} \mathbf{z}^\top) \left\{ \operatorname{MD}^2_{\boldsymbol{\mu}, \boldsymbol{\Sigma}}(\mathbf{x}) - \operatorname{trace}(\mathbf{P}_S \mathbf{z} \mathbf{z}^\top) \right\}.$$

Let $\theta$ denote the angle between $\mathbf{z}$ and $\mathbf{P}_S \mathbf{z}$ (its projection onto the $K$-eigenvector subspace) then $\cos(\theta) = \operatorname{trace}(\mathbf{P}_S \mathbf{z} \mathbf{z}^\top) / \{ \operatorname{MD}_{\boldsymbol{\mu}, \boldsymbol{\Sigma}}(\mathbf{x}) \| \mathbf{P}_S \mathbf{z} \| \}$ which then gives $\operatorname{trace}(\mathbf{P}_S \mathbf{z} \mathbf{z}^\top) = \operatorname{MD}^2_{\boldsymbol{\mu}, \boldsymbol{\Sigma}}(\mathbf{x}) \cos^2(\theta)$ (since $\| \mathbf{P}_S \mathbf{z} \| = \sqrt{\operatorname{trace}(\mathbf{P}_S \mathbf{z})}$) so that

$$\tilde{\rho}_S(C_0, F_{\boldsymbol{\mu}, \boldsymbol{\Sigma}}; \mathbf{x}) = \frac{\lambda_1 \lambda_p \operatorname{MD}^4_{\boldsymbol{\mu}, \boldsymbol{\Sigma}}(\mathbf{x})}{K(\lambda_1 - \lambda_p)^2} \cos^2(\theta) \{ 1 - \cos^2(\theta) \}.$$

In Figure 1 we plot $\tilde{\rho}_S(C_0, F_{\boldsymbol{\mu}, \boldsymbol{\Sigma}}; \mathbf{x})$ for Example 1 and the corresponding measures for the RMCD and $S$-estimators as described above. As can be seen in plot (a), the classical estimator can be highly influenced by $\mathbf{x}$, in particular those $\mathbf{x}$ with a large $\operatorname{MD}_{\boldsymbol{\mu}, \boldsymbol{\Sigma}}(\mathbf{x})$. However, the angle of $\mathbf{x}$ from its projection onto the subspace spanned by $[\boldsymbol{\eta}_1, \ldots, \boldsymbol{\eta}_K]$ also plays an important role. Regardless of the magnitude of $\operatorname{MD}_{\boldsymbol{\mu}, \boldsymbol{\Sigma}}(\mathbf{x})$, $\mathbf{x}$ has zero influence when $\mathbf{z} \in \operatorname{Span}\{\boldsymbol{\eta}_j : j \in S\}$ or $\mathbf{z} \in \operatorname{Span}\{\boldsymbol{\eta}_r : r \in S'\}$. This is also the case for the RMCD and $S$-estimator though the downweighting of observations with a large $\operatorname{MD}_{\boldsymbol{\mu}, \boldsymbol{\Sigma}}(\mathbf{x})$ results in reduced influence as is seen in plots (b) and (c). Let $q_\xi$ be the $\xi \times 100\%$ percentile for the $\chi_p^2$ distribution such that $P(\chi_p^2 \leq q_\xi) = \xi$. For the RMCD estimator there are discontinuities at $\operatorname{MD}_{\boldsymbol{\mu}, \boldsymbol{\Sigma}}(\mathbf{x}) = q_\alpha$ and $\operatorname{MD}_{\boldsymbol{\mu}, \boldsymbol{\Sigma}}(\mathbf{x}) = q_\delta$ corresponding to the rejection of observations in the initial weighting and re-weighting steps respectively. The estimator is particularly sensitive at these points since small changes in $\mathbf{x}$ can effect a large change in influence. From an influence perspective the $S$-estimator is preferred since it does not suffer from comparatively high influence and the smooth weighting function utilized results in smooth changes in influence with respect to changes in $\operatorname{MD}_{\boldsymbol{\mu}, \boldsymbol{\Sigma}}(\mathbf{x})$.

### *3.2. Dimension reduction methods*

In the regression setting, consider a univariate response variable $Y$ and $p$-dimensional predictor vector $\boldsymbol{X}$. If there exists a $p \times K$ matrix $\mathcal{B}$ such that $Y \perp\!\!\!\perp \boldsymbol{X} | \mathcal{B}^\top \boldsymbol{X}$ then, when $K < p$, dimension reduction can be achieved without loss of information by replacing the $p$-dimensional $\boldsymbol{X}$ with the $K$-dimensional $\mathcal{B}^\top \boldsymbol{X}$. For more information see, for e.g., [7]. Let $\mathcal{S}$ denote the column space of $\mathcal{B}$. Under appropriate conditions, methods such as Sliced Inverse Regression (SIR, [18]), Sliced Average Variance Estimates (SAVE, [6]) and Principal Hessian Directions (PHD, [19]) seek a basis for $\mathcal{S}$. The methods are based on an eigen-analysis of a $p \times p$ symmetric matrix so that, provided the influence function for the symmetric matrix estimator is known, the influence measure resulting from Theorem 1 is applicable in this setting.

Influence functions for versions of SIR, SAVE and PHD that return an orthonormal basis for $\mathcal{B}$ have been considered where influence on the directions of the basis is carried out with respect to Bénasséni's measure (see [23; 24] and [25]). As an example we will consider PHD since the method requires less introduction. Assuming $\boldsymbol{X} \sim N_p(\boldsymbol{\mu}, \boldsymbol{\Sigma})$, [19] showed that eigenvectors corresponding to non-zero eigenvalues of the the average Hessian matrix given as

$$\overline{\mathbf{H}}_{\mathbf{x}} = \boldsymbol{\Sigma}^{-1} E\left[\{Y - E(Y)\}\{\boldsymbol{X} - \boldsymbol{\mu}\}\{\boldsymbol{X} - \boldsymbol{\mu}\}^\top\right]\boldsymbol{\Sigma}^{-1}$$

are elements of $\mathcal{S}$. For more on PHD, including alternative versions, see [19] or [8].

In this regression context, the contamination distribution becomes

$$F_{\mathbf{x}}(\varepsilon) = (1 - \varepsilon)F_{\boldsymbol{\mu}, \boldsymbol{\Sigma}} + \varepsilon \Delta_{y, \mathbf{x}}$$

where $\Delta_{y, \mathbf{x}}$ has all of its probability mass at $(y, \mathbf{x}) \in \mathbb{R}^{p+1}$ to allow for contamination in both the response and predictor spaces. We shall assume that $F_{\boldsymbol{\mu}, \boldsymbol{\Sigma}} = N_p(\boldsymbol{\mu}, \boldsymbol{\Sigma})$ and $\text{rank}(\overline{\mathbf{H}}_{\mathbf{x}}) = K$ so that PHD is capable of finding a complete basis for $\mathcal{S}$. Let $H$ denote the functional for the usual estimator of $\overline{\mathbf{H}}_{\mathbf{x}}$. From [25] the influence function for $H$ at $F_{\boldsymbol{\mu}, \boldsymbol{\Sigma}}$ as

$$\begin{aligned}
\text{IF}(H, F_{\boldsymbol{\mu}, \boldsymbol{\Sigma}}; y, \mathbf{x}) &= \{y - E(Y)\}\left\{\mathbf{w}\mathbf{w}^\top - \boldsymbol{\Sigma}^{-1}\right\} - \mathbf{w}\left\{\mathbf{w}^\top \boldsymbol{\Sigma}\overline{\mathbf{H}}_{\mathbf{x}} + \mathbf{b}_{\text{OLS}}^\top\right\} \\
&\quad - \left\{\overline{\mathbf{H}}_{\mathbf{x}}\boldsymbol{\Sigma}\mathbf{w} + \mathbf{b}_{\text{OLS}}\right\}\mathbf{w}^\top - \overline{\mathbf{H}}_{\mathbf{x}}
\end{aligned} \tag{8}$$

where $\mathbf{w} = \boldsymbol{\Sigma}^{-1}(\mathbf{x} - \boldsymbol{\mu})$, $\mathbf{b}_{\text{OLS}}$ is the OLS slope vector at $F_{\boldsymbol{\mu}, \boldsymbol{\Sigma}}$ and $E(Y)$ is the expected value of $Y$ at $F_{\boldsymbol{\mu}, \boldsymbol{\Sigma}}$. We are interested in the basis estimator for $\mathcal{S}$ so we set $S = \{1, \ldots, K\}$ where each $\lambda_i$ ($i \in S$) is a non-zero eigenvalue with corresponding eigenvector $\boldsymbol{\eta}_i \in \mathcal{S}$. We also have (see, for e.g., [4; 5] or [22]) $\mathbf{b}_{\text{OLS}} \in \mathcal{S}$ so that, and since $\overline{\mathbf{H}}_{\mathbf{x}}\boldsymbol{\eta}_r = \mathbf{0}$ ($r \in S'$), from (5),

$$\begin{aligned}
\tilde{\rho}_S(H, F_{\boldsymbol{\mu}, \boldsymbol{\Sigma}}; y, \mathbf{x}) &= \frac{1}{K}\sum_{j \in S}\frac{1}{\lambda_j^2}\sum_{r \in S'}\Big[\{y - E(Y)\}\left(w_j w_r - \boldsymbol{\eta}_j^\top \boldsymbol{\Sigma}^{-1}\boldsymbol{\eta}_r\right) \\
&\quad - \left(\lambda_j \boldsymbol{\eta}_j^\top \boldsymbol{\Sigma}\mathbf{w} + \boldsymbol{\eta}_j^\top \mathbf{b}_{\text{OLS}}\right)w_r\Big]^2
\end{aligned} \tag{9}$$

where $w_m = \boldsymbol{\eta}_m^\top \mathbf{w}$.

The unboundedness of $\tilde{\rho}_S(H, F_{\boldsymbol{\mu}, \boldsymbol{\Sigma}}; y, \mathbf{x})$ is clearly evident suggesting that some observational types may be highly destructive to the estimator. In the original PHD paper by Li [19], it was noted that the method can be highly sensitive to outlying observations in the response space. The model itself imposes distributional restrictions on the predictor (i.e. normality) meaning that outliers in the predictor space may be more formally identified and subsequently treated. However, this is not the case with the response where outlying observations may still contain important regression information. However, $\tilde{\rho}_S(H, F_{\boldsymbol{\mu}, \boldsymbol{\Sigma}}; y, \mathbf{x})$ increases without bound as $y$ is moved further from $E(Y)$ suggesting that such observational types can be highly influential.

It is also interesting to note the types of observations that are not influential. For example, suppose that $\boldsymbol{\mu} = \mathbf{0}$ and $\boldsymbol{\Sigma} = \mathbf{0}$. Then, from (9), we have that $\tilde{\rho}_S(H, F_{\boldsymbol{\mu}, \boldsymbol{\Sigma}}; y, \mathbf{x}) = 0$ if $\mathbf{x}$ is an element of $\mathcal{S}$. That is, even extreme outliers may not be influential.

## 4. Sample versions

In practice it is common to consider the effect of highly influential observations on sample estimates. A limitation, however, is in the difficulty and inefficiency of locating such observations in large data sets. In this section we will consider sample based versions of the influence measure for the detection of influential observations.

Let $\mathbf{x}_1, \ldots, \mathbf{x}_n$ denote a random sample of size $n$ where each $\mathbf{x}_i \in \mathbb{R}^p$ and let $F_n$ denote the empirical distribution of this data. Similarly, suppose that $F_{n,(i)}$ denotes the empirical distribution of the data without the $i$th observation so that

$$F_{n,(i)} = \left(1 + \frac{1}{n-1}\right) F_n - \frac{1}{n-1} \Delta_{\mathbf{x}_i}.$$

Throughout, all reference will be to data of this form though in the regression setting one would need to consider observational pairs consisting of a predictor and response. Sample based versions of the influence function for the $i$th observation have been employed (see, for e.g., [9]) where the contaminant is $\mathbf{x}_i$ and $\epsilon$ is replaced with $-1/(n-1)$. For $\widehat{\mathbf{P}}$ and $\widehat{\mathbf{P}}_{S,(i)}$ denoting the projection matrix estimates with an without the $i$th observation, the true sample version of $\rho_S(W, F_{\boldsymbol{\mu}, \boldsymbol{\Sigma}}; \mathbf{x})$ is then

$$r_S(W, F_n; \mathbf{x}_i) = (n-1)^2 \left[1 - \frac{1}{K} \text{trace}(\widehat{\mathbf{P}}_S \widehat{\mathbf{P}}_{S,(i)})\right]. \tag{10}$$

Computation of $r_S(W, F_n; \mathbf{x}_i)$ for all $n$ observations requires estimation at each of $F_n, F_{n,(1)}, \ldots, F_{n,(n)}$. Such a process can be extremely inefficient when $n$ is large or $p$ is large (or both) due to the time it can take to carry out an eigen-analysis $n+1$ times. An approximation to $r_S(W, F_n; \mathbf{x}_i)$ may be computed by replacing unknown population parameters in $\tilde{\rho}_S(W, F_{\boldsymbol{\mu}, \boldsymbol{\Sigma}}; \mathbf{x})$ with their re-spective estimates at $F_n$. This approximate influence value is, from (5),

$$\tilde{r}_S(W, F_n; \mathbf{x}_i) = \frac{1}{K} \sum_{j \in S} \sum_{r \in S'} \frac{\{\widehat{\boldsymbol{\nu}}_j^\top \text{EIF}(W, F_n; \mathbf{x}_i) \widehat{\boldsymbol{\nu}}_r\}^2}{(\widehat{\kappa}_j - \widehat{\kappa}_r)^2} \tag{11}$$

where $\text{EIF}(W, F_n; \mathbf{x}_i)$ is the empirical influence function consisting of estimates at $F_n$ in place of population parameter values in $\text{IF}(W, F_{\boldsymbol{\mu}, \boldsymbol{\Sigma}}; \mathbf{x})$.

When $\text{IF}(W, F_{\boldsymbol{\mu}, \boldsymbol{\Sigma}}; \mathbf{x})$ exists in a closed form, i.e. in terms of $\mathbf{x}$ and population parameters only, then $\tilde{r}_S(W, F_n; \mathbf{x}_i)$ may be calculated for each observation after just one eigen-analysis at $F_n$. In the next section we will highlight the usefulness of this approximation in the context of computation time.

## 5. Sample principal components of the classical covariance matrix estimator: A microarray application

In this section we consider the colon tumor microarray data set [1]. For simplicity we consider the first $K$ estimated eigenvectors of the sample covariance matrix estimate $\mathbf{S}$. Computation of $r_S(W, F_n; \mathbf{x}_i)$ for even just a few observations can be extremely inefficient for high-dimensional data sets. If $r_S(W, F_n; \mathbf{x}_i)$ is to be calculated for all observations then such an analysis may become extremely onerous. We will therefore consider efficiently approximating $r_S(W, F_n; \mathbf{x}_i)$ with $\tilde{r}_S(W, F_n; \mathbf{x}_i)$. All results were obtained using R version 2.5.1 and the R function `eigen` for the eigen-analysis which utilizes LAPACK routines (see [2]) for computation. An Intel Pentium D CPU 3.60GHz with 1.99GB of RAM was used for the analysis.

The colon tumor microarray data set consists of gene expression measurements for 2000 genes corresponding to $n = 62$ samples. Each sample is either classified as being a 'normal tissue' sample or a 'tumor tissue' sample. Of the 40 individuals in the study, each has an associated 'tumor tissue' sample and 22 of the individuals also have a 'normal tissue' sample. We consider the normalized data where each sample is standardized to have mean 0 and standard deviation 1. Although this is a subset of a larger data set consisting of 6500 genes, most statistical research has concentrated on just the 2000 genes. We chose this data since it is often used in discriminant analysis but classical methods are not immediately applicable due to the singularity of $\mathbf{S}$. As such methods such as PCA may be used to initially reduce the dimension of the predictor space. For more information regarding this data set see [1].

From (6),

$$\tilde{r}_S(C_0, F_n; \mathbf{x}_i) = \frac{1}{K} \sum_{j \in S} \sum_{r \in S'} \frac{y_{ji}^2 y_{ri}^2}{(\widehat{\lambda}_j - \widehat{\lambda}_r)^2}. \tag{12}$$

where $y_{ji} = \widehat{\boldsymbol{\eta}}_j^\top (\mathbf{x}_i - \overline{\mathbf{x}})$, $\overline{\mathbf{x}}$ is the sample mean and $\widehat{\boldsymbol{\eta}}_1, \dots, \widehat{\boldsymbol{\eta}}_p$ are the sample eigenvectors corresponding to the sample eigenvalues $\widehat{\lambda}_1 \geq \widehat{\lambda}_2 \geq \cdots \geq \widehat{\lambda}_p$ of $\mathbf{S}$. Potential efficiency problems may still exist when the number of loop repetitions is large. A total of $K \times (p - K)$ iterations are required for the computation

TABLE 1

*Computation time in seconds for computation of all 62 $r_S(W, F_n; \mathbf{x}_i)$'s ($T_r$), $\tilde{r}_S(W, F_n; \mathbf{x}_i)$'s ($T_{\tilde{r}}$) and $\tilde{r}_S^*(W, F_n; \mathbf{x}_i)$'s ($T_{\tilde{r}}^*$) for the tumor data with $S = \{1, \ldots, K\}$. The spearman rank correlation between the $r_S(W, F_n; \mathbf{x}_i)$'s and $\tilde{r}_S^*(W, F_n; \mathbf{x}_i)$'s, $SR_S(r, \tilde{r})$, is also reported.*

| $S$ | $\{1\}$ | $\{1,2\}$ | $\{1,2,3\}$ | $\{1,\ldots,4\}$ | $\{1,\ldots,5\}$ |
|---|---|---|---|---|---|
| $T_r$ | 2953.24 | 2951.47 | 2950.50 | 2952.33 | 3021.73 |
| $T_{\tilde{r}}$ | 30.26 | 31.56 | 32.78 | 33.91 | 36.27 |
| $T_{\tilde{r}}^*$ | 29.41 | 29.37 | 29.50 | 29.52 | 29.56 |
| $SR_S(r,\tilde{r})$ | 0.995 | 0.993 | 0.975 | 0.929 | 0.962 |
| $S$ | $\{1,\ldots,6\}$ | $\{1,\ldots,7\}$ | $\{1,\ldots,8\}$ | $\{1,\ldots,9\}$ | $\{1,\ldots,10\}$ |
| $T_r$ | 3044.42 | 2953.89 | 2957.11 | 2955.84 | 2959.18 |
| $T_{\tilde{r}}$ | 36.11 | 37.34 | 39.08 | 40.13 | 41.60 |
| $T_{\tilde{r}}^*$ | 29.56 | 29.60 | 29.70 | 29.67 | 29.70 |
| $SR_S(r,\tilde{r})$ | 0.963 | 0.958 | 0.954 | 0.958 | 0.967 |
| $S$ | $\{1,\ldots,11\}$ | $\{1,\ldots,12\}$ | $\{1,\ldots,13\}$ | $\{1,\ldots,14\}$ | $\{1,\ldots,15\}$ |
| $T_r$ | 2957.14 | 2957.09 | 2962.31 | 2959.51 | 2958.94 |
| $T_{\tilde{r}}$ | 41.41 | 43.01 | 44.74 | 46.14 | 47.33 |
| $T_{\tilde{r}}^*$ | 29.78 | 29.74 | 29.75 | 29.78 | 29.77 |
| $SR_S(r,\tilde{r})$ | 0.958 | 0.960 | 0.940 | 0.904 | 0.913 |
| $S$ | $\{1,\ldots,16\}$ | $\{1,\ldots,17\}$ | $\{1,\ldots,18\}$ | $\{1,\ldots,19\}$ | $\{1,\ldots,20\}$ |
| $T_r$ | 2959.51 | 2961.66 | 2960.95 | 2961.72 | 2962.49 |
| $T_{\tilde{r}}$ | 49.16 | 48.81 | 50.53 | 53.01 | 53.53 |
| $T_{\tilde{r}}^*$ | 29.81 | 29.86 | 29.84 | 29.88 | 29.92 |
| $SR_S(r,\tilde{r})$ | 0.875 | 0.747 | 0.814 | 0.777 | 0.656 |

of a single $r_S(C_0, F_n; \mathbf{x}_i)$ so that the total number of iterations required for the computation of $r_S(C_0, F_n; \mathbf{x}_1), \ldots, r_S(C_0, F_n; \mathbf{x}_n)$ is $n \times K \times (p - K)$. For example, if we consider $n = 62$, $p = 2000$ and let $K = 10$, then the total number of iterations amounts to 1,239,380. However, this can be greatly reduced when $p >> n$ by noting that, since $\text{rank}(\mathbf{S}) \leq n - 1$ giving $\hat{\boldsymbol{\eta}}_k^\top \mathbf{S} \hat{\boldsymbol{\eta}}_k = 0$ for $k > n - 1$ then $y_{ki} = 0$ for all $i = 1, \ldots, n$ when $k = n, \ldots, p$. Hence

$$\tilde{r}_S^*(C_0, F_n; \mathbf{x}_i) = \tilde{r}_S(C_0, F_n; \mathbf{x}_i) = \frac{1}{K} \sum_{j \in S} \sum_{\substack{r \leq n-1 \\ r \notin S}} \frac{y_{ji}^2 y_{ri}^2}{(\hat{\lambda}_j - \hat{\lambda}_r)^2} \tag{13}$$

which requires just $K(n - 1 - K)$ iterations. Again for $n = 62$ and $p = 2000$ the total number of iterations required for a choice of $K = 10$ is only 31620; just 2.55% of the iterations required for (12).

In Table 1 we provide the time in seconds taken to compute all $r_S(W, F_n; \mathbf{x}_i)$'s and the approximations using $\tilde{r}_S(W, F_n; \mathbf{x}_i)$ and $\tilde{r}_S^*(W, F_n; \mathbf{x}_i)$. To highlight how well the approximation can be used to detect influential observations we also included the Spearman rank correlations between the $r_S(W, F_n; \mathbf{x}_i)$'s and $\tilde{r}_S^*(W, F_n; \mathbf{x}_i)$'s. It is immediately evident that much time can be saved when using the approximations. For example, for $K = 2$ the true computation cost is 2951.24 seconds compared to just 29.37 seconds (or around 1% of the time) using the $\tilde{r}_S^*(W, F_n; \mathbf{x}_i)$'s. The high Spearman rank correlation of 0.993 also indicates that $\tilde{r}_S^*(W, F_n; \mathbf{x}_i)$ is an excellent indicator of influential observations.

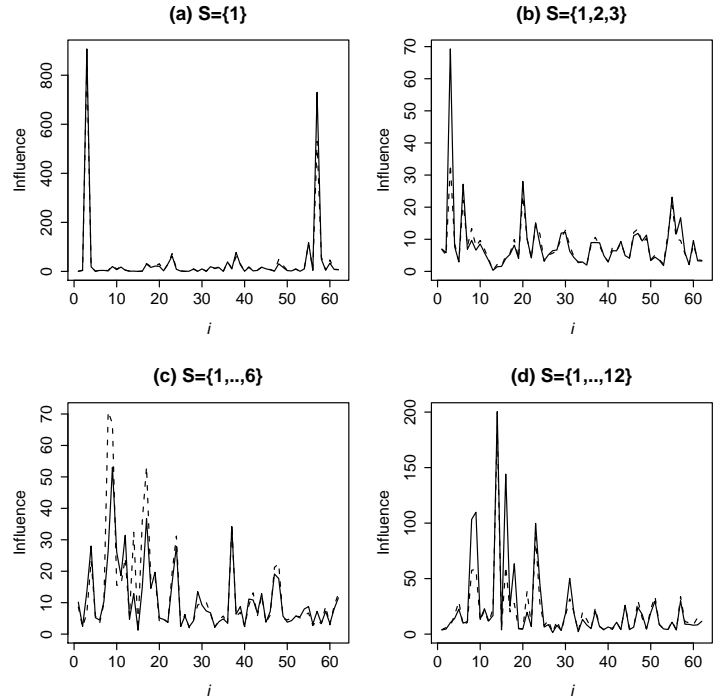

FIG 2. *Values of $r_S(W, F_n; \mathbf{x}_i)$'s (solid line) and $\tilde{r}_S^*(W, F_n; \mathbf{x}_i)$'s (dashed line) for the tumor data with (a) $S = \{1\}$, (b) $S = \{1, 2, 3\}$ (c) $S = \{1, \ldots, 6\}$ and (d) $S = \{1, \ldots, 12\}$.*

High correlations over 0.9 are maintained for all choices of $K$ up until and including $K = 15$. It is also clear that calculating the approximation based on the fewer loop iterations given in (13) is also beneficial with all computations coming in under 30 seconds compared to computation times approaching 60 seconds for larger choices of $K$ when (12) was utilized.

Further evidence of how close the approximation is to the true sample influence measure is provided in Figure 2. Here we plot the $r_S(W, F_n; \mathbf{x}_i)$'s versus the $\tilde{r}_S^*(W, F_n; \mathbf{x}_i)$'s for all observations and some choices of $K$. For $K = 1$ and 3 there is little difference between the true and approximate values. For $K = 6$ and 12, although there are obvious differences for some observations the approximation is still highly successful in highlighting influential observations.

## 6. Discussion

In this paper we considered sensitivity analysis of subsets of eigenvectors based on perturbation of the GCD and RV measures. Examples were provided to show how this analysis compliments existing influence studies in Section 3.

In Section 5 we highlighted how an approximate sample version of the measure can be used to efficiently detect influential observations in practice. The data set considered consisted of 2000 measurement variables for 62 individuals so that a leave-one-out sensitivity analysis requiring repeated principal component estimation was computationally inefficient. However, the approximate version provided an excellent approximation to the true sample measure when the subset of eigenvectors was not too large. This approximate version was much less time consuming to compute, therefore offering a useful means to assess influence for large data sets.

## Appendix A: Proof of Theorem 1

Using well-known eigen perturbation theory (see, for e.g., [26]), it is straightforward to show that $\mathbf{P}_S(\varepsilon)$ can be represented by the convergent power series

$$\mathbf{P}_S(\varepsilon) = \mathbf{P}_S + \varepsilon\mathbf{P}_1 + \frac{1}{2}\varepsilon^2\mathbf{P}_2 + O(\varepsilon^3). \tag{14}$$

For $S$ of the form $\{1, \ldots, K\}$, [31] give the form of $\mathbf{P}_1$ and $\mathbf{P}_2$. It is however, a simple generalization to use $S \subset \{1, \ldots, p\}$ and $S'$ in place of $\{1, \ldots, K\}$ and $\{K+1, \ldots, p\}$ respectively and simply let $K$ equal the number of elements in $S$.

Note that, since $(\mathbf{I} - \mathbf{P}_S)$ is a projection matrix,

$$
\begin{aligned}
K - \operatorname{trace}\left[\mathbf{P}_S\mathbf{P}_S(\varepsilon)\right] &= \operatorname{trace}\left[(\mathbf{I} - \mathbf{P}_S)\mathbf{P}_S(\varepsilon)\right] \\
&= \operatorname{trace}\left[(\mathbf{I} - \mathbf{P}_S)\mathbf{P}_S(\varepsilon)(\mathbf{I} - \mathbf{P}_S)\right] \\
&= \varepsilon\operatorname{trace}\left[(\mathbf{I} - \mathbf{P}_S)\left(\mathbf{P}_1 + \frac{1}{2}\varepsilon\mathbf{P}_2\right)(\mathbf{I} - \mathbf{P}_S)\right] + O(\varepsilon^3)
\end{aligned}
$$

from (14) and since $(\mathbf{I} - \mathbf{P}_S)\mathbf{P}_S = \mathbf{0}$.

The proof is complete by noting that, from [31], $(\mathbf{I} - \mathbf{P}_S)\mathbf{P}_1(\mathbf{I} - \mathbf{P}_S) = \mathbf{0}$ and

$$
\begin{aligned}
\operatorname{trace}\left[(\mathbf{I} - \mathbf{P}_S)\mathbf{P}_2(\mathbf{I} - \mathbf{P}_S)\right] &= \operatorname{trace}\left[2\sum_{r \in S'}\sum_{t \in S'}\sum_{j \in S}\frac{\boldsymbol{\nu}_j^\top\mathbf{W}_1\boldsymbol{\nu}_t}{\kappa_j - \kappa_t}\frac{\boldsymbol{\nu}_j^\top\mathbf{W}_1\boldsymbol{\nu}_r}{\kappa_j - \kappa_r}\boldsymbol{\nu}_r\boldsymbol{\nu}_t^\top\right] \\
&= 2\sum_{r \in S'}\sum_{j \in S}\frac{\left(\boldsymbol{\nu}_j^\top\mathbf{W}_1\boldsymbol{\nu}_r\right)^2}{\left(\kappa_j - \kappa_r\right)^2}
\end{aligned}
$$

since $\boldsymbol{\nu}_r^\top\boldsymbol{\nu}_t = 0$ $(r \neq t)$ and $\boldsymbol{\nu}_r^\top\boldsymbol{\nu}_r = 1$.